\documentclass{elsart}               % The use of LaTeX2e is preferred.

\usepackage{graphicx}          % Include this line if your
\begin{document}

\begin{frontmatter}

%\title{Optimal Design of Adiabatic Shortcuts for an Expanding Potential Well\thanksref{footnoteinfo}} % Title, preferably not more

\title{Optimal Shortcuts to Adiabaticity for a Quantum Piston\thanksref{footnoteinfo}}

\thanks[footnoteinfo]{The material of this paper has not been presented at any conference.}

\author[1]{Dionisis Stefanatos}
\ead{dionisis@post.harvard.edu}
\address{3 Omirou St., Sami, Kefalonia 28080, Greece}
\thanks[1]{Tel.: +30-697-4364682; Fax: +30-26740-23266.}

                             % keyword list or with the
                                          % help of the Automatica
                                          % keyword wizard

\begin{abstract}                          % Abstract of not more than 200 words.

In this paper we use optimal control to design minimum-time adiabatic-like paths for the expansion of a quantum piston. Under realistic experimental constraints, we calculate the minimum expansion time and compare it with that obtained from a state of the art inverse engineering method. We use this result to rederive the known upper bound for the cooling rate of a refrigerator, which provides a quantitative description for the unattainability of absolute zero, the third law of thermodynamics. We finally point out the relation of the present work to the fast adiabatic-like expansion of an accordion optical lattice, a system which can be used to magnify the initial quantum state (quantum microscope).

\end{abstract}

\end{frontmatter}

\section{Introduction}

Quantum science promises computational power unattainable by any classical computer \cite{Nielsen00}, as well as unprecedented precision measurement of various physical phenomena \cite{Wiseman09}. At the heart of these important applications lies the problem to accurately control and manipulate the states of quantum systems \cite{Dong10}. For many cases of interest, this control is achieved using adiabatic processes, where the system parameters are changed slowly from the initial to the desired final value \cite{Farhi01}. If the change is slow enough, the system follows the instantaneous eigenvalues and eigenstates of the time-dependent Hamiltonian and obtains the desired final state with a good approximation. The inherent drawback of adiabatic processes is that they require long times which may render them impractical, since most of the systems are not isolated but are exposed to undesirable interactions with the surrounding environment (decoherence) that lead to dissipation \cite{Breuer07}.

Several methods have been proposed to speed up adiabatic quantum dynamics. Their common characteristic is that they prepare a similar final state as the adiabatic process at a given final time (which in principle can be made arbitrarily short), without necessarily following the adiabatic path at each moment. The corresponding adiabatic-like trajectories are successfully characterized as shortcuts to adiabaticity. In the method of counter-diabatic control \cite{Demir08}, the applied electromagnetic field restores the adiabatic dynamics of the system by suppressing diabatic effects as they are generated. In the method of transitionless quantum driving \cite{Berry09}, an appropriate auxiliary time-dependent interaction is added such that the augmented system arrives in finite time at a similar quantum state with the adiabatic trajectory of the unperturbed system. In a closely related method, the trajectory of the shortcut is picked first and then the time-dependent interaction generating the corresponding evolution is inversely engineered \cite{Chen10}, using the theory of Lewis-Riesenfeld invariants \cite{Lewis69}. The above methods have been tested experimentally \cite{Schaff11,Bason12} and proven quite robust to various types of noise \cite{Choi12}. They also share another interesting feature: they do not specify a unique shortcut but rather provide entire families of them \cite{Chen11}. This ample freedom can be exploited using optimal control methods to find the shortcuts which minimize relevant physical criteria, like time, under realistic experimental constraints \cite{Stefanatos11}.

In this paper we use optimal control theory to design the shortest adiabatic-like shortcut for a potential well with a moving boundary, in the presence of restrictions suggested by the experimental setup. This system provides a model for a \emph{quantum piston} and has some interesting applications. Note that the shortcuts to adiabaticity for this quantum piston have been studied in \cite{AdolfoBox}, as well as their potential applications to the control of many-body quantum dynamics. Here we mostly concentrate on the control aspects of the problem, but we also highlight some nice applications. In Section \ref{sec:shortcuts} we summarize the important points of the above work which are essential for the current analysis and also formulate the time-optimal control problem for the expansion of the piston. We solve this problem in Section \ref{sec:shortest} and obtain the minimum expansion time as a function of the expansion factor of the piston. In Section \ref{sec:results} we compare our optimal results with those obtained using the inverse engineering method. Using the fact that the minimum expansion time has a logarithmic dependence on the expansion factor for large values of the latter, we reobtain the cooling rate of a refrigerator \cite{Hoffmann_EPL11}, a result which quantifies the unattainability of absolute zero implied by the third law of thermodynamics \cite{Rezek11}. Finally, we discuss how the control problem studied here is related to the fast adiabatic-like expansion of an optical lattice with dynamically variable spacing \cite{Yuce12}. Section \ref{sec:conclusion} concludes the paper.

\section{Shortcuts to Adiabaticity for a Quantum Piston}
\label{sec:shortcuts}

\subsection{Fast Adiabatic-Like Expansion}

Consider the potential
\begin{equation}
\label{well}
V_p(x,t)=\left\{\begin{array}{cc} 0, & 0<x<a(t)\\\infty, & \mbox{otherwise}\end{array}\right.
\end{equation}
which describes the infinite square well with a moving boundary at $x=a(t)$. The evolution of the wavefunction $\psi(x,t)$ of a particle trapped in this potential is given by the following Schr\"{o}dinger equation
\begin{equation}
\label{Schrodinger}
i\hbar\frac{\partial\psi}{\partial t}=\left[-\frac
{\hbar^{2}}{2m}\frac{\partial^{2}}{\partial x^{2}}+V_p(x,t)\right] \psi,
\end{equation}
where $m$ is the particle mass and $\hbar$ is Planck's constant. The wavefunction $\psi$ is square-integrable on the interval $[o,a]$, with $|\psi(x,t)|^2dx$ expressing the probability to find the particle between $x$ and $x+dx$ at time $t$. The above physical system can serve as a model for a quantum piston.

We are interested in the expansion of the piston from $a(0)=a_0>0$ to $a(T)=a_T>a_0$, at some final time $t=T$. If the expansion is slow enough (adiabatic) then, according to the adiabatic theorem, the solution of (\ref{Schrodinger}) is given approximately by
\begin{equation}
\label{adiabatic}
\psi(x,t)\approx\sum_{n=1}^\infty c_n\exp{\left[i\phi_n(t)\right]}\Psi_n(x,t),
\end{equation}
where
\begin{equation}
\label{eigen}
\Psi_n(x,t)=\sqrt{\frac{2}{a(t)}}\sin{\left[\frac{n\pi x}{a(t)}\right]}
\end{equation}
are the instantaneous eigenstates of the right hand side in (\ref{Schrodinger}) \cite{Landau}, $\phi_n$ are the so-called adiabatic phases, and $c_n$ are constant coefficients determined by the initial condition. The slower is the expansion, the better is the approximation in (\ref{adiabatic}). For fast expansion the adiabatic approximation breaks down.

A method has been proposed recently to accelerate adiabatic quantum dynamics, according to which an auxiliary potential is designed such that the system acquires a similar quantum state with the reference adiabatic path (\ref{adiabatic}) in arbitrarily short time $T$ \cite{Berry09}. For the case of an expanding piston, it has been shown in \cite{AdolfoBox} that if the following potential is added in (\ref{Schrodinger})
%\begin{equation}
%\label{auxiliary}
%V_{\mbox{aux}}(x,t)=\frac{1}{2}k(t)x^2,\quad k(t)=-m\frac{\ddot{a}(t)}{a(t)},
%\end{equation}
\begin{eqnarray}
\label{auxiliary}V_a(x,t) &=& \frac{1}{2}k(t)x^2,\\
\label{k}k(t) &=& -m\frac{\ddot{a}(t)}{a(t)},
\end{eqnarray}
then the \emph{exact} solution of the resulting Schr\"{o}dinger equation can be expressed as
\begin{eqnarray}
\label{exact}
\psi(x,t)&=&\sum_{n=1}^\infty c_n\exp{\left[-i\int_0^tE_n(t')dt'/\hbar\right]}\times\\
&&\exp{\left[i\frac{m\dot{a}(t)}{2\hbar a(t)}x^2\right]}\Psi_n(x,t),\nonumber
\end{eqnarray}
where
\begin{equation}
\label{energies}
E_n(t)=\frac{n^2\pi^2\hbar^2}{2ma^2(t)}
\end{equation}
are the instantaneous energy eigenvalues corresponding to the eigenstates (\ref{eigen}) \cite{Landau}.
The boundary conditions
\begin{equation}
\label{boundary}
a(0)=a_0,a(T)=a_T,\quad \dot{a}(0)=\dot{a}(T)=0,
%a(t)=\left\{\begin{array}{cc} a_0, & t=0\\a_T, & t=T\end{array}\right.,\quad \dot{a}(t)=\left\{\begin{array}{cc} 0, & t=0\\0, & t=T\end{array}\right.,
\end{equation}
ensure the expansion of the piston from $a_0$ to $a_T$ in the interval $0<t<T$, and that the exact solution (\ref{exact}) is similar to the adiabatic one (\ref{adiabatic}) at $t=0$ and $t=T$. The additional conditions
\begin{equation}
\label{kboundary}
%k(t)=\left\{\begin{array}{cc} 0, & t\leq 0\\0, & t\geq T\end{array}\right.
k(0)=k(T)=0
\end{equation}
guarantee that the auxiliary potential is active only within $0<t<T$.

\subsection{Inverse Engineering Approach}

The auxiliary potential can be designed using an elegant inverse engineering method \cite{Chen10}. The final time $T$ is fixed and the condition (\ref{kboundary}) is translated through (\ref{k}) to
\begin{equation}
\label{ddotaboundary}
%\ddot{a}(t)=\left\{\begin{array}{cc} 0, & t\leq 0\\0, & t\geq T\end{array}\right..
\ddot{a}(0)=\ddot{a}(T)=0.
\end{equation}
A polynomial ansatz satisfying (\ref{boundary}) and (\ref{ddotaboundary}) has been found in \cite{AdolfoBox}, $a(\tau)/a_0=1+(\gamma-1)\tau^3(6\tau^2-15\tau+10)$, where $\tau=t/T$ and $\gamma=a_T/a_0$ is the expansion factor. Then, the stiffness $k(t)$ of the auxiliary potential can be determined from (\ref{k})
\begin{equation}
\label{kshortcut}
k(\tau)=-\frac{m}{T^2}\frac{60(\gamma-1)\tau(2\tau^2-3\tau+1)}{1+(\gamma-1)\tau^3(6\tau^2-15\tau+10)},
\end{equation}
where again $\tau=t/T$. Note that there is no mathematical limitation on the size of $T$, which can be chosen arbitrarily small in theory.

\subsection{Optimal Control Approach}

In practise, there are always experimental constraints, for example
\begin{equation}
\label{kconstraint}
-k_0\leq k(t)\leq k_0,
\end{equation}
which restrict $T$ to some finite value. In such cases, finding the shortest adiabatic-like path can be expressed as an optimal control problem. If we set
\begin{equation}
x_{1}=\frac{a(t)}{a_0},\quad x_{2}=\sqrt{\frac{m}{k_0}}\frac{\dot{a}(t)}{a(t)},\quad u(t)=\frac{k(t)}{k_0},\nonumber
\end{equation}
and normalize time according to $t_{\mathrm{new}}=t_{\mathrm{old}}/T_0$, where $T_0=\sqrt{m/k_0}$, we
obtain the following system, equivalent
to equation (\ref{k})
\begin{eqnarray}
\label{system1}\dot{x}_{1} &=& x_{2},\\
\label{system2}\dot{x}_{2} &=& -ux_{1}.
\end{eqnarray}
The minimum time adiabatic shortcut, under the constraint (\ref{kconstraint}), can be found by solving the following optimal control problem:
\begin{prob}
\label{problem}
Find $u(t)$ with $-1\leq u(t)\leq 1$ such that starting from $(x_1(0),x_2(0))=(1,0)$, the above system with $x_1(t)>0$ reaches the final point $(x_1(T),x_2(T))=(\gamma,0)$, $\gamma>1$, in minimum time $T$.
\end{prob}

The boundary conditions on $x_1$ and $x_2$ correspond to those for $a$ and $\dot{a}$ from (\ref{boundary}), while the additional constraint $x_1(t)>0$ corresponds to the natural requirement $a(t)>0$. In the following section we solve the above problem on the interval $0<t<T$. In order to much the boundary conditions $u(0)=u(T)=0$, corresponding to (\ref{kboundary}), the optimal control may be complemented with instantaneous jumps at the initial and final times which do not affect the cost (time), see Fig. \ref{fig:controls} in Section \ref{sec:results}. This approach is similar to that used in our recent work \cite{Stefanatos11,Stefanatos12}, as well as in \cite{Hoffmann_EPL11}. Note that the problem of finding the control $u(t)$ with $0<u_{\mathrm{min}}\leq u\leq u_{\mathrm{max}}$ which drives in minimum time the system (\ref{system1}) and (\ref{system2}) (the parametric oscillator) from the point $(x_1,x_2)$ to the ellipse $x_2^2+u_Tx_1^2=2E_T$, with specified oscillator energy $E_T$ for stiffness $u_T$, has been studied thoroughly in \cite{Andresen11}. In Problem \ref{problem} the control is allowed to take negative values, while the target is a point and not a curve of the state space.

%%%%%%%%%%%%%%%%%%%%%%%%%%%%%%%%%%%%%%%%%%%%%%%%%%%%%%%%%%%
\section{The Shortest Shortcut to Adiabaticity}
\label{sec:shortest}

The system described by (\ref{system1}) and (\ref{system2}) can be expressed
in compact form as
\begin{equation}
\label{affine}
\dot{x}=f(x)+ug(x),
\end{equation}
where the vector fields are given by
\begin{equation}
f=\left[\begin{array}[c]{c}x_2\\0\end{array}\right],\quad
g=\left[\begin{array}[c]{c}0\\-x_1\end{array}\right]
\end{equation}
and $x\in\mathcal{D}=\{(x_1,x_2)\in R^2:x_1>0\}$, $u\in
U=[-1,1]$. Admissible controls are Lebesgue measurable functions that
take values in the control set $U$. Given an admissible control $u$ defined
over an interval $[0,T]$, the solution $x$ of the system (\ref{affine})
corresponding to the control $u$ is called the corresponding trajectory and we
call the pair $(x,u)$ a controlled trajectory.

For a constant $\lambda_{0}$ and a row vector $\lambda\in(R^2)^\ast$ the control
Hamiltonian for system (\ref{affine}) is defined as $H=H(\lambda_{0},\lambda,x,u)=\lambda_{0}+\lambda [f(x)+ug(x)]$. Pontryagin's Maximum Principle \cite{Pontryagin}
provides the following necessary optimality conditions:
\begin{thm}[Maximum principle]
\label{prop:max_principle}
Let $(x_{\ast}(t),u_{\ast}(t))$
be a time-optimal controlled trajectory that transfers the initial condition
$x(0)=x_0$ of system (\ref{affine}) into the terminal state $x(T)=x_T$. Then it is a necessary
condition for optimality that there exists a constant $\lambda_{0}\leq0$ and
nonzero, absolutely continuous row vector function $\lambda(t)$ such that:
\begin{enumerate}
\item $\lambda$ satisfies the adjoint equation $\dot{\lambda}=-\partial H/\partial x$.
\item For $\,0\leq t\leq T$ the function $u\mapsto H(\lambda_{0},\lambda(t),x_{\ast}(t),u)$ attains its maximum\ over the control set $U$ at
$u=u_{\ast}(t)$.
\item $H(\lambda_{0},\lambda(t),x_{\ast}(t),u_{\ast}(t))\equiv0$.
\end{enumerate}
\end{thm}

In the following we use maximum principle to solve Problem \ref{problem}.

\begin{defn}
We denote the vector fields corresponding to the constant bang controls
$u=-1$ and $u=1$ by $X=f-g$ and $Y=f+g$, respectively, and call
the corresponding trajectories $X$- and $Y$-trajectories. A concatenation of an
$X$-trajectory followed by a $Y$-trajectory is denoted by $XY$ while the
concatenation in the inverse order is denoted by $YX$.
\end{defn}

\begin{thm}[Optimal solution]
\label{optimal_thm}
The optimal trajectory for Problem \ref{problem} has the one-switching form $XY$. The optimal control is
\begin{equation}
\label{optimal}
u(t)=\left\{\begin{array}{cc} -1, & 0<t<T_X\\1, & T_X<t<T_X+T_Y\end{array}\right.,
\end{equation}
where
\begin{eqnarray}
\label{Tx}T_X &=& \sinh^{-1}\left(\sqrt{\frac{\gamma^2-1}{2}}\right),\\
\label{Ty}T_Y &=& \sin^{-1}\left(\frac{1}{\gamma}\sqrt{\frac{\gamma^2-1}{2}}\right).
\end{eqnarray}
The minimum expansion time is
\begin{equation}
\label{min_time}
T=T_X+T_Y
\end{equation}
\end{thm}
\begin{pf*}{Proof.}
We show first that the optimal control is bang-bang, i.e., alternates between the boundary values $u=\pm1$ of the control set. For the system (\ref{system1}), (\ref{system2}) we have
\begin{equation}
\label{hamiltonian}
H(\lambda_{0},\lambda,x,u)=\lambda_0+\lambda_1x_2-\lambda_{2}x_{1}u,
\end{equation}
and thus
\begin{eqnarray}
\label{la1}\dot{\lambda}_{1} &=& u\lambda_{2},\\
\label{la2}\dot{\lambda}_{2} &=& -\lambda_{1}.
\end{eqnarray}
Observe that $H$ is a linear function of the bounded control variable $u$.
The coefficient at $u$ in $H$ is $-\lambda_{2}x_{1}$ and, since $x_{1}>0$, its
sign is determined by $\Phi=-\lambda_{2}$, the so-called \emph{switching
function}. According to the maximum principle, point 2 above, the optimal
control is given by $u=\mbox{sign}\,\Phi$, if $\Phi\neq0$. The
maximum principle provides a priori no information about the control at times
$t$ when the switching function $\Phi$ vanishes. Now observe that whenever $\Phi(t)=-\lambda_{2}(t)=0$ at some
time $t$, then $\dot{\Phi}(t)=-\dot{\lambda}_{2}(t)=\lambda_{1}(t)\neq 0$ since the maximum principle requires that it is always $\lambda=(\lambda_1,\lambda_2)\neq 0$. Hence, when $\Phi(t)=0$ it is also
$\dot{\Phi}(t)\neq0$ and there is a switch between the control boundary values at time $t$. The optimal trajectory consists of a concatenation of $X$- and $Y$-trajectories.

We now move to narrow the candidate sequences for optimality. We show first that the concatenation $XYX$ cannot be part of the optimal trajectory. Without loss of generality assume that the $XY$ switch takes place at $t=0$, so $\lambda_2(0)=0$, while it is also $x_1(0)>0$. Along the $Y$-trajectory it is $u=1$, so from (\ref{la1}) and (\ref{la2}) we find $\lambda_2(t)=-\lambda_1(0)\sin t$. The subsequent $YX$ switch should take place at $t=\pi$, where the switching function $\Phi=-\lambda_2$ becomes zero again. But observe that for $u=1$ the state equations (\ref{system1}) and (\ref{system2}) correspond to a rotation around the origin with period $2\pi$, so at $t=\pi$ the state $x$ has been rotated by half circle and consequently $x_1(\pi)<0$, which is forbidden. We next show similarly that the sequence $YXY$ cannot also be part of the optimal trajectory. We assume that the $YX$ switch takes place at $t=0$, thus $\lambda_2(0)=0$. Along the $X$-trajectory it is $u=-1$, so from (\ref{la1}) and (\ref{la2}) we obtain $\lambda_2(t)=-\lambda_1(0)\sinh t$. Observe that $\Phi=-\lambda_2\neq 0$ for $t>0$, and the subsequent $XY$ switch is not allowed.

We conclude that the only candidates for optimality left are the $XY$ and $YX$ trajectories. Furthermore, it is not hard to see that only the former corresponds to expansion (final $\gamma>1$), while the latter corresponds to compression and is rejected. The equations of the $X$-trajectory starting from $(1,0)$ and of the $Y$-trajectory ending at $(\gamma,0)$ can be found from (\ref{system1}) and (\ref{system2}) for $u=\mp 1$ and they are
\begin{eqnarray}
\label{X}x_1^2-x_2^2 &=& 1,\\
\label{Y}x_1^2+x_2^2 &=& \gamma^2.
\end{eqnarray}
The switching point satisfies both equations and it is $(\sqrt{\gamma^2+1}/\sqrt{2},\sqrt{\gamma^2-1}/\sqrt{2})$, thus using (\ref{system1}) and (\ref{system2}) the times $T_X$ and $T_Y$ spent on each segment of the optimal trajectory can be easily derived as in (\ref{Tx}) and (\ref{Ty}).\qed
\end{pf*}

%%%%%%%%%%%%%%%%%%%%%%%%%%%%%%%%%%%%%%%%%%%%%%%%%%%%%%%%%%%
\section{Results and Discussion}
\label{sec:results}

\begin{figure}[t]
\begin{center}
\includegraphics[height=5.5cm]{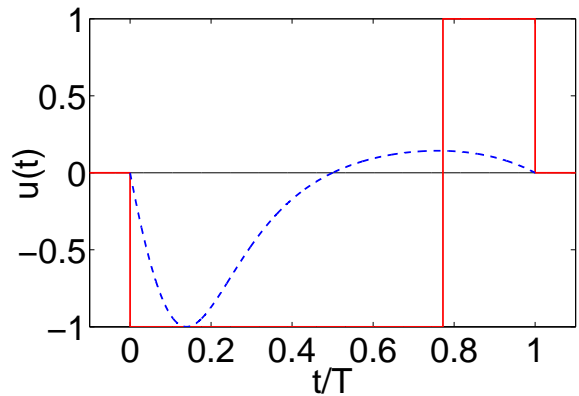}    % The printed column
\caption{Optimal control (solid line, $T=3.4295T_0$) and inverse engineering control (dashed line, $T=6.2511T_0$), under the constraint $-1\leq u(t)\leq 1$ and for the expansion of the piston by a factor of $\gamma=10$. Observe that both controls are active only within $0<t<T$.}                         % width is 8.4 cm.
\label{fig:controls}                           % Size the figures
\end{center}                               % accordingly.
\end{figure}

In Fig. \ref{fig:controls} we plot the time-optimal control $u(t)$ from (\ref{optimal}) (solid line, $T=3.4295T_0$), as well as the control $k(t)/k_0$ obtained from the inverse engineering method (\ref{kshortcut}) and with minimum duration $T$ under the constraint (\ref{kconstraint}) (dashed line, $T=6.2511T_0$), both for the final expansion factor $a_T/a_0=\gamma=10$. Observe that the control derived from the inverse engineering method is actually limited by the lower bound $u=-1$.

\begin{figure}[t]
\begin{center}
\includegraphics[height=5.5cm]{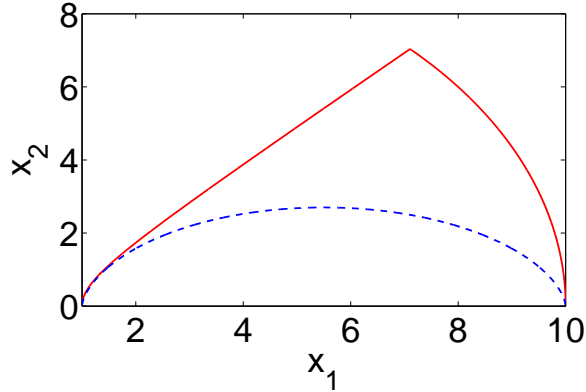}    % The printed column
\caption{Corresponding trajectories for the control inputs of Fig. \ref{fig:controls}.}                         % width is 8.4 cm.
\label{fig:trajectories}                           % Size the figures
\end{center}                               % accordingly.
\end{figure}

In Fig. \ref{fig:trajectories} we depict the corresponding trajectories. Note that the system equations (\ref{system1}) and (\ref{system2}) can be interpreted as describing the one-dimensional Newtonian motion of a unit-mass particle, with position coordinate $x_{1}$ and velocity $x_{2}$. The acceleration (force) acting on the particle is $-ux_{1}$. As we can observe from Fig. \ref{fig:trajectories}, in both trajectories the particle traverses the same distance $x_1$ but along the time-optimal trajectory its speed $\dot{x}_1=x_2$ is always higher.

\begin{figure}[t]
\begin{center}
\includegraphics[height=5.5cm]{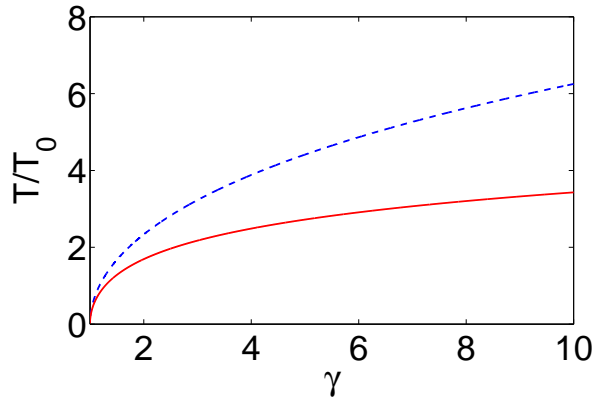}    % The printed column
\caption{Expansion time $T$ (in units of $T_0=\sqrt{m/k_0}$) as a function of the expansion factor $\gamma$, for the optimal (solid line) and the inverse engineering (dashed line) strategies.}                         % width is 8.4 cm.
\label{fig:times}                           % Size the figures
\end{center}                               % accordingly.
\end{figure}

In Fig. \ref{fig:times} we plot the expansion time $T$ in units of $T_0=\sqrt{m/k_0}$ as a function of the expansion factor $\gamma$ for both control policies. This graph provides the speed limits of the adiabatic-like expansion under condition (\ref{kconstraint}). Even if the optimal bang-bang controls are not experimentally exactly realizable, knowledge of the time-optimal solutions is a useful guide for the design of more realistic controls. For example, the abrupt changes of the optimal control shown in Fig. \ref{fig:controls} can be approximated by ramps of finite duration \cite{Stefanatos10,Hoffmann_EPL11}. In general, more complicated constraints on the control or even the state can be incorporated in the current formalism using a powerful numerical optimization method based on pseudospectral approximations \cite{Stefanatos10}.%$T_r$, resulting to a suboptimal trapezoidal control. We have performed numerical simulations which show that if the rise time is a reasonable fraction of $T_Y^\infty=\pi/4$ (the limiting value of $T_Y$ as $\gamma\rightarrow\infty$ in (\ref{Ty})), for example $T_r=T_Y^\infty/3=\pi/12$, then the corresponding expansion time is quite robust and close to the lower bound plotted in Fig. \ref{fig:times}.

The above results have an interesting thermodynamic application. The quantum piston using noninteracting particles as the working medium and executing the Otto cycle provides a model for a refrigerator, similar to that considered in \cite{Rezek09}. Following the procedure described in this work, the cooling rate $R$ of the refrigerator can be calculated, as the temperature of the cold reservoir $\tau_c$ approaches absolute zero. This rate is defined as $R=Q/T$, the ratio of the heat $Q$ extracted from the cold reservoir on each cycle to the duration $T$ of the cycle. We find $Q$ first. When the working medium is in contact with the hot reservoir of temperature $\tau_h$, its internal energy is $\tau_h/2$ (ideal gas with one degree of freedom). The next step in the Otto cycle is the adiabatic expansion of the piston from $a_0$ to $a_T$. At the end of this process the populations of the energy levels are preserved (\ref{adiabatic}), while the energies are reduced by a factor of $\gamma^2=a_T^2/a_0^2$ (\ref{energies}), thus the internal energy of the working medium becomes $\tau_h/2\gamma^2$. After the expansion the working medium is brought in contact with the cold reservoir and its internal energy is raised to $\tau_c/2$. The heat extracted from the cold reservoir is $Q=\tau_c/2-\tau_h/2\gamma^2$. A necessary condition for operation of the refrigerator is $Q>0\Rightarrow\gamma>\sqrt{\tau_h/\tau_c}$. For $\tau_c\rightarrow 0$ this implies also $\gamma\rightarrow\infty$, and in this limit the duration of the cycle is dominated by the duration of the adiabatic expansion. From (\ref{Tx}) and (\ref{Ty}) we find that for the fastest adiabatic-like expansion it is $T=T_X+T_Y\rightarrow\ln\gamma$ for $\gamma\rightarrow\infty$, thus the cooling rate is restricted as $R=Q/T<-\tau_c/\ln\tau_c$ for $\tau_c\rightarrow 0$, and the result of \cite{Hoffmann_EPL11} is reobtained. As $\tau_c\rightarrow 0$ the cooling rate approaches zero faster, a manifestation of the unattainability of the absolute zero (third thermodynamics law) \cite{Rezek11}.

We finally mention another possible application of the present work. The optimal control problem solved in the previous section is directly related to the fast adiabatic-like expansion of an optical lattice with dynamically variable spacing \cite{Yuce12}. The connection is that the stiffness $\omega^2$ of the external harmonic potential, necessary to keep the atoms trapped during the expansion, is related to the lattice scale parameter $\Lambda$ by the relation $\omega^2=-\ddot{\Lambda}/\Lambda$, which is similar to (\ref{k}). As a consequence, the present analysis applies also to this context. Note that such accordion lattices are useful since the final lattice spacing can be made large enough to be resolved experimentally, for example by imaging of the atoms at individual sites. The final quantum state is a scaled-up version of the initial state due to the adiabatic-like evolution, thus the optical lattice acts like a quantum dynamical microscope \cite{AdolfoMicroscope}.

%%%%%%%%%%%%%%%%%%%%%%%%%%%%%%%%%%%%%%%%%%%%%%%%%%%%%%%%%%%
\section{Conclusion}
\label{sec:conclusion}

In this paper we formulated and solved the problem of minimum-time adiabatic-like expansion for a quantum piston, in the presence of experimental constraints. As a result, we obtained the speed limit for this fast quantum driving, and used it to rederive an interesting result related to the third law of thermodynamics. We also highlighted the possible application of the present work to the adiabatic-like expansion of an optical lattice.

%\section{The argument}
%\begin{equation} \label{e1}
%{{\partial F}\over {\partial t}} =
%D{{\partial^2 F}\over {\partial x^2}}.
%\end{equation}

%\begin{thm}
%The square of the length of the hypotenuse of a right triangle equals the sum of the squares
%of the lengths of the other two sides.
%\end{thm}
%\section{Epilogue}
%A word or two to conclude, and this even includes some inline
%maths:  $R(x,t)\sim t^{-\beta}g(x/t^\alpha)\exp(-|x|/t^\alpha)$.

%\begin{ack}                               % Place acknowledgements
%Partially supported by the Roman Senate.  % here.
%\end{ack}

%\bibliographystyle{plain}        % Include this if you use bibtex
%\bibliography{autosam}           % and a bib file to produce the
                                 % bibliography (preferred). The
                                 % correct style is generated by
                                 % Elsevier at the time of printing.

%\parskip 0pt

%\appendix
%\section{A summary of Latin grammar}    % Each appendix must have a short title.
%\section{Some Latin vocabulary}         % Sections and subsections are supported
                                        % in the appendices.
\end{document}